\def\autori{P.\ DALL'AGLIO}
\def\titolo{Asymptotic behaviour of solutions of fourth order
Dirichlet problems} 




\font\sixrm=cmr6
\newcount\tagno \tagno=0                        
\newcount\thmno \thmno=0                        
\newcount\bibno \bibno=0                        
\newcount\chapno\chapno=0                       
\newcount\verno            
\newif\ifproofmode
\proofmodefalse
\newif\ifwanted
\wantedfalse
\newif\ifindexed
\indexedfalse

\def\ifundefined#1{\expandafter\ifx\csname+#1\endcsname\relax}

\def\Wanted#1{\ifundefined{#1} \wantedtrue 
\immediate\write0{Wanted 
#1 
\the\chapno.\the\thmno}\fi}

\def\Increase#1{{\global\advance#1 by 1}}

\def\Assign#1#2{\immediate
\write1{\noexpand\expandafter\noexpand\def
 \noexpand\csname+#1\endcsname{#2}}\relax
 \global\expandafter\edef\csname+#1\endcsname{#2}}

\def\pAssign#1#2{\write1{\noexpand\expandafter\noexpand\def
 \noexpand\csname+#1\endcsname{#2}}}

\def\lPut#1{\ifproofmode\llap{\hbox{\sixrm #1\ \ \ }}\fi}
\def\rPut#1{\ifproofmode$^{\hbox{\sixrm #1}}$\fi}



\def\chp#1{\global\tagno=0\global\thmno=0\Increase\chapno
\Assign{#1}
{\the\chapno}{\lPut{#1}\the\chapno}}


\def\thm#1{\Increase\thmno 
\Assign{#1}
{\the\chapno.\the\thmno}\the\chapno.\the\thmno\rPut{#1}}


\def\frm#1{\Increase\tagno
\Assign{#1}{\the\chapno.\the\tagno}\lPut{#1}
{\the\chapno.\the\tagno}}


\def\bib#1{\Increase\bibno
\Assign{#1}{\the\bibno}\lPut{#1}{\the\bibno}}


\def\pgp#1{\pAssign{#1/}{\the\pageno}}


\def\ix#1#2#3{\pAssign{#2}{\the\pageno}
\immediate\write#1{\noexpand\idxitem{#3}
{\noexpand\csname+#2\endcsname}}}

 
\def\rf#1{\Wanted{#1}\csname+#1\endcsname\relax\rPut {#1}}


\def\rfp#1{\Wanted{#1}\csname+#1/\endcsname\relax\rPut{#1}}

\input \jobname.auxi
\Increase\verno
\immediate\openout1=\jobname.auxi

\immediate\write1{\noexpand\verno=\the\verno}

\ifindexed
\immediate\openout2=\jobname.idx
\immediate\openout3=\jobname.sym 
\fi


\font\twelverm=cmr12
\font\twelvei=cmmi12
\font\twelvesy=cmsy10
\font\twelvebf=cmbx12
\font\twelvett=cmtt12
\font\twelveit=cmti12
\font\twelvesl=cmsl12

\font\ninerm=cmr9
\font\ninei=cmmi9
\font\ninesy=cmsy9
\font\ninebf=cmbx9
\font\ninett=cmtt9
\font\nineit=cmti9
\font\ninesl=cmsl9

\font\eightrm=cmr8
\font\eighti=cmmi8
\font\eightsy=cmsy8
\font\eightbf=cmbx8
\font\eighttt=cmtt8
\font\eightit=cmti8
\font\eightsl=cmsl8

\font\sixrm=cmr6
\font\sixi=cmmi6
\font\sixsy=cmsy6
\font\sixbf=cmbx6

\catcode`@=11 
\newskip\ttglue

\def\twelvepoint{\def\rm{\fam0\twelverm}
\textfont0=\twelverm  \scriptfont0=\ninerm  
\scriptscriptfont0=\sevenrm
\textfont1=\twelvei  \scriptfont1=\ninei  \scriptscriptfont1=\seveni
\textfont2=\twelvesy  \scriptfont2=\ninesy  
\scriptscriptfont2=\sevensy
\textfont3=\tenex  \scriptfont3=\tenex  \scriptscriptfont3=\tenex
\textfont\itfam=\twelveit  \def\it{\fam\itfam\twelveit}%
\textfont\slfam=\twelvesl  \def\sl{\fam\slfam\twelvesl}%
\textfont\ttfam=\twelvett  \def\tt{\fam\ttfam\twelvett}%
\textfont\bffam=\twelvebf  \scriptfont\bffam=\ninebf
\scriptscriptfont\bffam=\sevenbf  \def\bf{\fam\bffam\twelvebf}%
\tt  \ttglue=.5em plus.25em minus.15em
\normalbaselineskip=15pt
\setbox\strutbox=\hbox{\vrule height10pt depth5pt width0pt}%
\let\sc=\tenrm  \let\big=\twelvebig  \normalbaselines\rm}

\def\tenpoint{\def\rm{\fam0\tenrm}
\textfont0=\tenrm  \scriptfont0=\sevenrm  \scriptscriptfont0=\fiverm
\textfont1=\teni  \scriptfont1=\seveni  \scriptscriptfont1=\fivei
\textfont2=\tensy  \scriptfont2=\sevensy  \scriptscriptfont2=\fivesy
\textfont3=\tenex  \scriptfont3=\tenex  \scriptscriptfont3=\tenex
\textfont\itfam=\tenit  \def\it{\fam\itfam\tenit}%
\textfont\slfam=\tensl  \def\sl{\fam\slfam\tensl}%
\textfont\ttfam=\tentt  \def\tt{\fam\ttfam\tentt}%
\textfont\bffam=\tenbf  \scriptfont\bffam=\sevenbf
\scriptscriptfont\bffam=\fivebf  \def\bf{\fam\bffam\tenbf}%
\tt  \ttglue=.5em plus.25em minus.15em
\normalbaselineskip=12pt
\setbox\strutbox=\hbox{\vrule height8.5pt depth3.5pt width0pt}%
\let\sc=\eightrm  \let\big=\tenbig  \normalbaselines\rm}

\def\ninepoint{\def\rm{\fam0\ninerm}
\textfont0=\ninerm  \scriptfont0=\sixrm  \scriptscriptfont0=\fiverm
\textfont1=\ninei  \scriptfont1=\sixi  \scriptscriptfont1=\fivei
\textfont2=\ninesy  \scriptfont2=\sixsy  \scriptscriptfont2=\fivesy
\textfont3=\tenex  \scriptfont3=\tenex  \scriptscriptfont3=\tenex
\textfont\itfam=\nineit  \def\it{\fam\itfam\nineit}%
\textfont\slfam=\ninesl  \def\sl{\fam\slfam\ninesl}%
\textfont\ttfam=\ninett  \def\tt{\fam\ttfam\ninett}%
\textfont\bffam=\ninebf  \scriptfont\bffam=\sixbf
\scriptscriptfont\bffam=\fivebf  \def\bf{\fam\bffam\ninebf}%
\tt  \ttglue=.5em plus.25em minus.15em
\normalbaselineskip=11pt
\setbox\strutbox=\hbox{\vrule height8pt depth3pt width0pt}%
\let\sc=\sevenrm  \let\big=\ninebig  \normalbaselines\rm}

\def\eightpoint{\def\rm{\fam0\eightrm}
\textfont0=\eightrm  \scriptfont0=\sixrm  \scriptscriptfont0=\fiverm
\textfont1=\eighti  \scriptfont1=\sixi  \scriptscriptfont1=\fivei
\textfont2=\eightsy  \scriptfont2=\sixsy  \scriptscriptfont2=\fivesy
\textfont3=\tenex  \scriptfont3=\tenex  \scriptscriptfont3=\tenex
\textfont\itfam=\eightit  \def\it{\fam\itfam\eightit}%
\textfont\slfam=\eightsl  \def\sl{\fam\slfam\eightsl}%
\textfont\ttfam=\eighttt  \def\tt{\fam\ttfam\eighttt}%
\textfont\bffam=\eightbf  \scriptfont\bffam=\sixbf
\scriptscriptfont\bffam=\fivebf  \def\bf{\fam\bffam\eightbf}%
\tt  \ttglue=.5em plus.25em minus.15em
\normalbaselineskip=9pt
\setbox\strutbox=\hbox{\vrule height7pt depth2pt width0pt}%
\let\sc=\sixrm  \let\big=\eightbig  \normalbaselines\rm}

\def\twelvebig#1{{\hbox{$\textfont0=\twelverm\textfont2=\twelvesy
	\left#1\vbox to10pt{}\right.\n@space$}}}
\def\tenbig#1{{\hbox{$\left#1\vbox to8.5pt{}\right.\n@space$}}}
\def\ninebig#1{{\hbox{$\textfont0=\tenrm\textfont2=\tensy
	\left#1\vbox to7.25pt{}\right.\n@space$}}}
\def\eightbig#1{{\hbox{$\textfont0=\ninerm\textfont2=\ninesy
	\left#1\vbox to6.5pt{}\right.\n@space$}}}
 
\def\displayliness#1{\null\,\vcenter{\openup1\jot \m@th
  \ialign{\strut\hfil$\displaystyle{##}$\hfil
    \crcr#1\crcr}}\,}
	       
\def\displaylinesno#1{\displ@y \tabskip=\centering
   \halign to\displaywidth{ \hfil$\@lign \displaystyle{##}$ \hfil
	\tabskip=\centering
     &\llap{$\@lign##$}\tabskip=0pt \crcr#1\crcr}}
		     
\def\ldisplaylinesno#1{\displ@y \tabskip=\centering
   \halign to\displaywidth{ \hfil$\@lign \displaystyle{##}$\hfil
	\tabskip=\centering
     &\kern-\displaywidth
     \rlap{$\@lign##$}\tabskip=\displaywidth \crcr#1\crcr}}

\catcode`@=12 

\def\parag#1#2{\goodbreak\bigskip\bigskip\noindent
		   {\bf #1.\ \ #2}
		   \nobreak\bigskip} 
\def\intro#1{\goodbreak\bigskip\bigskip\goodbreak\noindent
		   {\bf #1}\nobreak\bigskip\nobreak}
\long\def\th#1#2{\goodbreak\bigskip\noindent
		{\bf Theorem #1.\ \ \it #2}}
\long\def\lemma#1#2{\goodbreak\bigskip\noindent
		{\bf Lemma #1.\ \ \it #2}}
\long\def\prop#1#2{\goodbreak\bigskip\noindent
		  {\bf Proposition #1.\ \ \it #2}}

\long\def\rem#1#2{\goodbreak\bigskip\noindent
		 {\bf Remark #1.\ \ \rm #2}}
\long\def\ex#1#2{\goodbreak\bigskip\noindent
		 {\bf Example #1.\ \ \rm #2}}

\def\proof{\vskip.4cm\noindent{\it Proof.\ \ }}

\def\sqr#1#2{\vbox{
   \hrule height .#2pt 
   \hbox{\vrule width .#2pt height #1pt \kern #1pt 
      \vrule width .#2pt}
   \hrule height .#2pt }}
\def\square{\sqr74}

\def\endproof{{\unskip\nobreak\hfill \penalty50
\hskip2em\hbox{}\nobreak\hfill $\square$ \goodbreak
\parfillskip=0pt  \finalhyphendemerits=0}}

\mathchardef\emptyset="001F
\mathchardef\hyphen="002D


\def\wto{\rightharpoonup}


\def\rightheadline{\eightpoint\hfil\titolo
\hfil\tenrm\folio} 
\def\leftheadline{\tenrm\folio\hfil\eightpoint
\autori \hfil}
\def\zeroheadline{\hfill} 
%
\headline={\ifnum\pageno=0 \zeroheadline
\else\ifodd\pageno\rightheadline
\else\leftheadline\fi\fi}

\nopagenumbers
\magnification=1200
\baselineskip=15pt
\hfuzz=2pt
\parindent=2em
\mathsurround=1pt
\tolerance=1000

\pageno=0
\hsize 14truecm
\vsize 25truecm
\hoffset=0.8truecm
\voffset=-1.55truecm

\null
\vskip 2.8truecm
{\twelvepoint
\baselineskip=1.7\baselineskip
\centerline{\bf ASYMPTOTIC BEHAVIOUR OF SOLUTIONS}
\centerline{\bf OF FOURTH ORDER DIRICHLET PROBLEMS}
}
\vskip2truecm

\centerline{Paolo DALL'AGLIO}
\vfil

{\eightpoint
\baselineskip=1.2\baselineskip
\centerline{\bf Abstract}
\bigskip
\noindent 
The behaviour of solutions to fourth order problems is studied through 
the decomposition into a system of second order ones, which leads to 
relaxed formulations with the introduction of measure terms. This allows 
to solve a shape optimization problem for a simply supported thin 
plate.
\par } 
\vfil 
\vskip 1truecm 
\centerline {Ref. S.I.S.S.A. 
32/97/M (March 97)}
\vskip 1truecm 
\eject
%
%

%
%
\topskip=25pt 
\vsize 22.5truecm
\hsize 16.2truecm
\hoffset=0truecm
\voffset=0.5truecm


\def\R{\rm I\! R}
\def\Ldue#1{{\rm L}^2(#1)}
\def\Lduemu#1{{\rm L}^2_\mu(#1)}
\def\Hunozero{{\rm H}^1_0}
\def\Hduale{{\rm H}^{\hbox{\kern 1truept \rm -}\kern -1truept1}\kern 
	-1truept} 
\def\Mo{{\cal M}_0(\Omega)} 
\def\intl{\int\limits}
\def\intersA#1{{\rm H}^1_0(#1)\cap {\rm L}^2_{\mu_A}(#1)}
\def\intersB#1{{\rm H}^1_0(#1)\cap {\rm L}^2_{\mu_B}(#1)}
\def\dualita#1#2{\langle#1,#2\rangle} 
\def\Re{R\kern -1truept e}

\parag{\chp{intro}}{Introduction}

In this paper we study the asymptotic behaviour of solutions of fourth 
order elliptic problems on varying domains.

This has been widely studied in the past, in the case of second order 
elliptic operators (see for instance [\rf{CIO-MUR}], [\rf{DM-MOS}],
[\rf{DM-GAR}]). We will use such results decomposing fourth order 
differential equations into a system of second order ones.

Given a bounded open set $U$ in $\R^n$, $n\geq2$ and a function 
$f\in\Hduale(\Omega)$, the fourth order equation 
$$
\cases{	\Delta^2 u=f \hbox{ in }\Hduale(\Omega)\cr
	\Delta u\in\Hunozero(\Omega)\cr
	u\in\Hunozero(\Omega)\cr}\eqno(\frm{probl})
$$
is linked to the model for  the vertical displacement $u$ of an
thin plate, occupying a region $U$, 
simply supported on $\partial U$, subjected to a load $f$.
Simply supported means that the boundary is fixed, but that the plate is 
free to 
rotate around the tangent to $\partial U$. For the general treatment of
plate theory we refer to [\rf{LAN-LIF}], 
[\rf{DUV-LIO}], [\rf{CIA1}], [\rf{LAG-LIO}], [\rf{CIA2}].

In particular we want to study the asymptotic behaviour of solutions  
when the domain varies.
To this aim we will show (Proposition \rf{equiv}) that problem 
(\rf{probl}) is equivalent to the system of second order equations
$$
\cases{-\Delta u=v\hbox{ in }\Hduale(U)\cr
	u\in\Hunozero(U)\cr
	-\Delta v=f\hbox{ in }\Hduale(U)\cr
	v\in\Hunozero(U).\cr}\eqno(\frm{probl3})
$$
This problem can be handled with the theorems valid in the second order 
case (Theorems \rf{compactness} and \rf{density}), and it will be proved 
that if $u_n$ are the solutions of 
problems like (\rf{probl3}) on a sequence of subdomains $U_n$ of a given
bounded domain $\Omega$, then a
subsequence of $u_n$ converges weakly in $\Hunozero(\Omega)$ to a function 
$u_\mu$ solving $$
\cases{-\Delta u_\mu+\mu u_\mu=v\cr
	u_\mu\in\Hunozero(\Omega)\cap\Lduemu{\Omega}\cr
	-\Delta v+\mu v=f\cr
	v\in\Hunozero(\Omega)\cap\Lduemu{\Omega}.\cr}\eqno(\frm{pblim})
$$
where $\mu$ is a measure.

The study will be carried on for general fourth order elliptic operators 
with constant coefficients and no lower order terms, that can be splitted 
into two second order ones.

A motivation for the study of the asymptotic behaviour of solutions of
Dirichlet problems in varying domains without geometric assumptions on the
domains $U_n$ are the so-called shape 
optimization problems:
given a function $j:\Omega\times\R\to\R$ 
we consider the following problem
$$
\min_{U\in{\cal U}(\Omega)}\intl_\Omega j(x,u_{{}_U}(x))\, dx,\eqno(\frm{optim})
$$
where ${\cal U}(\Omega)$ is the family of all open subset of $\Omega$ and
$u_{{}_U}$ is the solution of the problem of type (\rf{probl}) in the 
set $U$.

In section \rf{optimization} it will be shown that, in general, problem 
(\rf{optim}) does not have a solution. Hence a relaxed optimization 
problem will be introduced, where the set over which we minimize is the 
set of functions $u_\mu$, where $\mu$ is a measure and $u_\mu$ solves 
the relaxed problem (\rf{pblim}). This set is the closure of 
$\left\{u_{{}_U}\,:\,U\in{\cal U}(\Omega)\right\}$ in $\rm L^2(\Omega)$. 
This problem will always have solution and its minimum will coincide with 
the infimum of integral in (\rf{optim}).

This shape optimization problem, in the second order case, was studied in 
[\rf{ATT}], [\rf{BUT-DM}], [\rf{CHI-DM}].

\parag{\chp{not}}{Notations and preliminary results}

Given an open subset $U$ of $\R^n$, $\Hunozero(U)$ is the usual 
Sobolev Space, $\Hduale(U)$ its dual, and $\dualita{\cdot}{\cdot}$ the 
duality pairing. If $U\subset\Omega$ and 
$u\in\Hunozero(U)$, then the function 
$$
\tilde u:=\cases{u &in $U$\cr
		0 &in $\Omega\setminus U$\cr}
$$
is in $\Hunozero(\Omega)$. From now on we will always denote, with the same 
symbol $u$, a function and its extension $\tilde u$.

In this paper we will deal with elliptic operators $L:{\cal 
D'}(\Omega)\rightarrow{\cal D'}(\Omega)$ of the form 
$$
Lu=\sum_{|\alpha|=m} 
c_{\alpha}\partial^\alpha u,
$$
where $\alpha$ and $\beta$ are multiindeces, and $m$ is the order of the 
operator that will be 2 or 4. In any case they will be without lower order 
term and with constant coefficients.

The operators will be assumed to be elliptic in the sense that
$$
\sum_{|\alpha|=m}c_{\alpha}\xi^\alpha\ge\gamma|\xi|^m,
\qquad\forall\,\xi\in\R^n,
$$
where $\gamma$ is a real positive constant. This, in our case, is 
the same as
$$
P(\xi)\neq0\quad\forall\,\xi\in\R^n\setminus\{0\},\eqno(\frm{ellipt})
$$
where $P$ is the polynomial $\displaystyle\sum_{|\alpha|=m}
c_{\alpha}\xi^\alpha$ associated to the operator $L$.

In this work, differential problems are always meant to be solved 
in the usual 
weak sense. This means, for instance, that, for $u\in \Hunozero(U)$ the 
expression
$$
-\Delta u=f\hbox{ in }\Hduale(U)
$$ 
is an equality of linear functionals
$$
\dualita{-\Delta u}{v}=\sum_i\intl_U\partial_i
u\partial_iv\,dx=\dualita{f}{v},
$$
for any $v\in\Hunozero(U)$.

As said above, the limit of a sequence of solutions of Dirichlet problems
is not, in general, 
the solution of a problem of the same kind, but is the solution of a 
problem where a measure term appears. To deal with these problems we 
need to recall some notions.

For the notion of capacity of a set $E\subset\Omega$, which we will 
indicate by  cap$(E)$, we refer to textbooks as [\rf{EVA-GAR}] or 
[\rf{HEI-KIL}]. We shall always identify a function 
$u\in\Hunozero(\Omega)$ with its quasi-continuous representative.

Now, let $\Mo$ be the set of Borel measures which are zero on the sets of 
zero capacity. 

For one such measure $\mu\in\Mo$, $\Lduemu{U}$ will be the space 
of functions such that
$$
\intl_U |u|^2\,d\mu\,<\,+\infty.
$$
Given a second order operator $A=\displaystyle{\sum_{ij}a_{ij}\partial_i 
\partial_j}$, with constant coefficients, 
for a function $u\in\Hunozero(U)\cap\Lduemu{U}$ to solve the equation
$$
Au +\mu u=f,
$$
will mean that
$$
\sum_{i=1}^n\intl_U a_{ij}\partial_j u \partial_i v\,dx+\intl_U uv\,d\mu = 
	\intl_U fv\,dx,\eqno(\frm{relax})
$$
for all test functions in $v\in\Hunozero(U)\cap\Lduemu{U}$.

It can be easily proved that the space $\Hunozero(U)\cap\Lduemu{U}$ is a 
Hilbert space whenever $\mu$ is in $\Mo$, and hence, by Lax-Milgram Lemma, 
we have existence and uniqueness \pgp{E!} of solutions for a problem of the 
form 
$$
\cases{Au+\mu u=f\cr
	u\in\Hunozero(U)\cap\Lduemu{U},\cr}
$$
for any linear elliptic second order operator $A$.

The decomposition of a fourth order problem in a system of two second
order equations allows us to study 
the asymptotic behaviour applying 
well known theorems for the second order case 
separately to each equation. The following  results, that 
can be found, for instance 
in [\rf{DM-MOS}], [\rf{DM-GAR}], are the key points of the theory.

\th{\thm{compactness}}{Let $A$ be a second order elliptic operator, as
described above. 
For every sequence $\{\mu_n\}$ in $\Mo$ there 
exists a subsequence $\mu_{n_k}$ such that, for every
sequence $\{g_n\}$ in 
$\Hduale(\Omega)$, strongly convergent to $g\in\Hduale(\Omega)$, we have
$$
z_{n_k}\wto z\hbox{ weakly in }\Hunozero(\Omega),
$$
where $z_{n_k}$ and $z$ solve
$$
\cases{Az_{n_k}+\mu_{n_k} z_{n_k}=g_{n_k}\cr
	z_{n_k}\in\Hunozero(\Omega)\cap 
	L^2_{\mu_{n_k}}(\Omega),\cr}\qquad
\cases{Az+\mu z=g\cr                                     
	u\in\Hunozero(\Omega)\cap\Lduemu{\Omega},\cr}
$$
respectively.}

\th{\thm{density}}{For any measure $\mu\in\Mo$ there exists a sequence 
$U_n$ of open subsets of $\Omega$ such that, for any sequence $\{g_n\}$ in 
$\Hduale(\Omega)$, strongly convergent to $g\in\Hduale(\Omega)$, 
the solutions $z_n$ of the problems
$$
\cases{Az_n=g_n \quad\hbox{ in  }\Hduale(U_n)\cr
	z_n\in\Hunozero(U_n)\cr}
$$
converge weakly in $\Hunozero(\Omega)$ to the solution $z$ of
$$
\cases{Az+\mu z=g\cr                                     
	u\in\Hunozero(\Omega)\cap\Lduemu{\Omega}.\cr}
$$
}

\parag{\chp{res}}{Decomposition of fourth order operators}

Fourth order elliptic problems  are studied mainly with two different 
kinds of boundary conditions. In the model case 
of the bi-laplacian, they correspond to two different physical
problems regarding, as said above, the displacement of a thin plate. 
Problems of the type
$$
\cases{	\Delta^2 u=f \hbox{ in }\Hduale(\Omega)\cr
	u\in {\rm H}^2_0(\Omega)\cr}
$$
correspond to having $u={\displaystyle{{\partial u}\over{\partial n}}}=0$
on $\partial\Omega$, that is to say that the plate is clamped along its 
boundary. The asymptotic behaviour of the solutions of such problems 
has been studied in [\rf{DM-PAD}]. 

In this work we deal with the second
kind of boundary conditions, as in (\rf{probl}), and we do it decomposing 
the problem in a system of two second order equations.

The decomposability of the fourth order operator 
$$
Lu=\sum_{|\alpha|=4} 
c_{\alpha}\partial^\alpha u
$$
can be seen through the associated polynomial 
$$
P(\xi)=\sum_{|\alpha|=4}c_{\alpha}\xi^\alpha.
$$

It is a simple algebraic fact that, if the polynomial can be splitted
into two second degree polynomials
$$
P(\xi)=Q(\xi)R(\xi),
$$
then other decompositions can be obtained only by exchanging the order 
or multiplying and dividing by constants.
Observe that, according to (\rf{ellipt}), if $P$ is elliptic, 
then so are $R$ and $Q$.

We remark here that such a decomposition can always be done in 
the two dimensional case. In higher dimensions this is not always
possible.

So assume that
$$
Lu=\sum_{|\alpha|=4}c_{\alpha}\partial ^\alpha \,=\,
\sum_{i,j=1}^n b_{ij}\partial _j \partial _i
\left(\sum_{k,l=1}^n a_{kl}\partial _k \partial _l\right)u=BAu.
$$

\prop{\thm{equiv}}{Let $U$ be a subset of $\Omega$, and 
$A={\displaystyle\sum_{i,j=1}^n a_{ij}\partial _j \partial _i}$, 
$B={\displaystyle\sum_{i,j}^n b_{ij}\partial _j \partial _i}$ be 
second order elliptic operators with constant coefficients. Let
$f\in\Hduale(U)$. The following three problems are equivalent:}
\medskip
$$
(i)\cases{BAu=f {\rm \ \ in\ } \Hduale(U)\cr
       Au\in\Hunozero(U)\cr
       u\in\Hunozero(U)\cr}
(ii)\cases{u\in\Hunozero(U)\ :\ Au\in\Ldue{U}\cr
	\displaystyle{\int_U AuB\varphi\,dx }= \langle f,\varphi\rangle\cr
	\forall\varphi\in\Hunozero(U)\ :\ B\varphi\in\Ldue{U}\cr}
(iii)\cases{Au=v {\rm\ \ in\ }\Hduale(U)\cr
	u\in\Hunozero(U)\cr
	Bv = f{\rm\ \ in\ }\Hduale(U)\cr
	v\in\Hunozero(U)\cr}\eqno(\frm{trepb})
$$

\proof
$(i)\Rightarrow (ii)$. The equality $\langle BAu,\varphi\rangle=\langle 
f,\varphi\rangle$ holds in 
particular for any $\varphi$ in $\Hunozero(U)$ such that 
$B\varphi\in\Ldue{U}$. Since $B$ is symmetric and $Au\in\Hunozero(U)$ we 
get 
$$
\intl_U AuB\varphi\,dx=\langle B\varphi,Au\rangle=\langle
BAu,\varphi\rangle=
\dualita{f}{\varphi}
$$

$(ii)\Rightarrow (iii)$. Let $v$ be solution of
$$
v\in\Hunozero(U),\quad
\dualita{Bv}{\varphi}=\dualita{f}{\varphi}\quad\forall\varphi\in\Hunozero(U).
$$
Then
$$
\dualita{B\varphi}{v}=\dualita{f}{\varphi}
$$
and if, in particular, $\varphi$ is such that $B\varphi\in\Ldue{U}$, we get
$$
\intl_U vB\varphi\,dx=\langle f,\varphi\rangle,
$$
which, subtracted to the equation in $(ii)$ gives
$$
\intl_U(Au-v)B\varphi\,dx=0.
$$
Observe now that, thanks to Lax-Milgram theorem, every  
function in $\Ldue{U}$ can 
be written as $B\varphi$, with $\varphi\in\Hunozero(U)$, so we obtain
$$
\intl_U(Au-v)z\,dx=0,\ \ \forall z\in\Ldue{U}.
$$
Taking $z=Au-v$, we get that $\Vert Au-v\Vert_{\Ldue{U}}^2=0$, hence
$Au=v$, and $u$ solves problem $(iii)$.

$(iii)\Rightarrow (ii)$. If $u$ is a solution of $(iii)$ then
$$
\dualita{Bv}{\varphi}=\dualita{f}{\varphi}\quad\forall\varphi\in\Hunozero(U).
$$
and if, in particular, $B\varphi\in\Ldue{U}$ we have
$$
\intl_U vB\varphi\,dx=\langle f,\varphi\rangle.
$$
Since now  $Au=v$ we have
$$
\intl_U AuB\varphi\,dx = \langle f,\varphi\rangle,
$$
for every $\varphi\in\Hunozero(U)$ with $B\varphi\in\Ldue{U}$, hence $u$ 
solves $(ii)$.

$(ii)\Rightarrow (i)$. We have already proved that from $(ii)$ it follows 
that $Au\in\Hunozero(U)$, hence $BAu\in\Hduale(U)$, and so
$$
\langle f,\varphi\rangle=\intl_U AuB\varphi\,dx=
\langle B\varphi,Au\rangle=\dualita{BAu}{\varphi};
$$
that is, $BAu=f$ in $\Hduale(U)$\endproof

\rem{\thm{esist}}{} The equivalence af these problems gives for free 
existence and uniqueness, since this is true for problem $(iii)$ 
thanks to the observation made in section \rf{not}. This is not obvious
for 
problem $(ii)$ because the spaces of solution and of test functions are 
different.

\rem{\thm{hduezero}}{}Notice that if the boundary of $U$ is regular, then 
by regularity theorems, $Au\in\Ldue{U}$ implies that $u$ belongs to 
${\rm H}^2(U)$. The same is true for the test functions $\varphi$. So the 
problem(\rf{trepb})$(ii)$ becomes:
$$
\cases{u\in\Hunozero(U)\cap {\rm H}^2(U)\cr
	\displaystyle{\int_U AuB\varphi\,dx }= \langle f,\varphi\rangle\cr
	\forall\varphi\in\Hunozero(U)\cap {\rm H}^2(U).\cr}
$$

\rem{\thm{diversi}}{} It is important to remark that the boundary 
conditions in the equivalent problems (\rf{trepb})$(i)(ii)(iii)$ depend 
on the choice of the decomposition $L=BA$.

If we have two decompositions, in the sense that
$$
L\varphi\,=\,B_1A_1\varphi\,=\,B_2A_2\varphi,\qquad\forall\varphi\in{\cal 
D}'(\Omega), 
$$
then the boundary conditions in (\rf{trepb})$(i)$ are different. Hence 
also the other problems $(ii)$ and $(iii)$ differ.
Had we chosen the Dirichlet boundary conditions, that is 
seeking $u\in {\rm H}^2_0(U)$, all decompositions would have given the same 
solution. But in this case what makes the difference is the 
boundary condition, as can be seen with the following example.

\ex{\thm{ex}}{} The operator 
$$
Lu\,:= u_{xxxx}+2u_{yyyy}+3u_{xxyy}
$$
can be written, for instance, as the product of
$$
Au=\Delta u\ \hbox{ and }\ Bu=\Delta u+u_{yy}.
$$
But the integral in equation (\rf{trepb})$(ii)$ will  be 
different according to which operator we will apply first:
$$
\displaylines{
\int Au Bv\,dx\,=\,\int( 
	u_{xx}v_{xx}+2u_{yy}v_{yy}+2u_{xx}v_{yy}+u_{yy}v_{xx})\,dx\, ,\cr
\int Bu Av\,dx\,=\,\int (
	u_{xx}v_{xx}+2u_{yy}v_{yy}+u_{xx}v_{yy}+2u_{yy}v_{xx})\,dx\, .\cr} 
$$

Computations show that these two integrals differ by a term on the boundary, 
which would vanish if functions where in ${\rm H}^2_0(U)$.

\parag{\chp{limits}}{The asymptotic behaviour}

We come now to examine problem (\rf{trepb})$(i)$ when we have 
a sequence of domains $U_n$ all contained in $\Omega$.
Let $f$ be in $\Hduale(\Omega)$ (this implies it is also in 
$\Hduale(U_n)$ for any $U_n$) and consider always 
functions of $\Hunozero(U_n)$ trivially extended to the whole of $\Omega$.
As proved in 
Proposition \rf{equiv}, we can study directly 
$$
\cases{Au_n=v_n {\rm\ \ in\ }\Hduale(U_n)\cr
	u_n\in\Hunozero(U_n)\cr
	Bv_n = f{\rm\ \ in\ }\Hduale(U_n)\cr
	v_n\in\Hunozero(U_n).\cr}\eqno(\frm{systemn})
$$
We first consider the second equation
$$
\cases{Bv_n = f{\rm\ \ in\ }\Hduale(U_n)\cr
	v_n\in\Hunozero(U_n).\cr}
$$
From Theorem \rf{compactness} we know that there exist a subsequence, which 
we still call $U_n$,  and a measure $\mu_B$, depending on the sequence 
$U_n$, on the operator $B$, but not on $f$, such that 
$$
v_n\wto v \ \ \hbox{ weakly in }\Hunozero(\Omega)
$$
and $v$ solves
$$
\cases{Bv+\mu_B v=f\cr
	v\in\intersB{\Omega},\cr}\eqno(\frm{probleB})
$$
in the sense specified in (\rf{relax}).

We can now apply Theorem \rf{compactness} to the problem in $u$
$$
\cases{Au_n = v_n{\rm\ \ in\ }\Hduale(U_n)\cr
	u_n\in\Hunozero(U_n),\cr}
$$
taking as $U_n$ only those in the subsequence obtained for $B$.
Again there exists a subsequence, which 
we still call $U_n$,  and a measure $\mu_A$, depending on the sequence 
$U_n$, on the operator $A$, but not on the sequence $v_n$, such that 
$$
u_n\wto u \hbox{ weakly in }\Hunozero(\Omega)
$$
and $u$ solves
$$
\cases{Au+\mu_A u=v\cr
	u\in\intersA{\Omega}.\cr}\eqno(\frm{probleA})
$$
This allows us to conclude that, if $u_n$ are the solutions of
system (\rf{systemn}) then, up to a subsequence,
$$
u_n\wto u \hbox{ weakly in }\Hunozero(\Omega)
$$
and $u$ is the solution of
$$
\cases{Au+\mu_A u=v\cr
	u\in\intersA{\Omega}\cr
	Bv+\mu_B v=f\cr
	v\in\intersB{\Omega}.\cr}\eqno(\frm{problelimi})
$$

\parag{\chp{equation}}{The single equation formulation}

The goal of this section is to write problem (\rf{problelimi}) as a single 
equation of the form
$$
\intl_{\Omega}(Au+\mu_A u)(B\varphi +\mu_B\varphi)\,dx=\langle
f,\varphi\rangle.
$$

Of course $Au+\mu_A u$ alone doesn't make sense, because it has to be 
understood in the sense explained in chapter \rf{not}.
What we will do now, is hence to define suitable function spaces, which 
will play the role of 
$\left\{z\in{\Hunozero(\Omega)}\,:\,Az\in\Ldue{\Omega}\right\}$ and 
$\left\{w\in{\Hunozero(\Omega)}\,:\,Bw\in\Ldue{\Omega}\right\}$ and 
operators $\Lambda_{\mu_A}$ and $\Lambda_{\mu_B}$ which give meaning
to the integral
$$
\intl_\Omega \Lambda_{\mu_A}u\ \Lambda_{\mu_B}\varphi\,dx.
$$
The construction of such operators will be done only relatively to 
operator $A$, being the one relative to $B$ perfectly analogous.

Consider the operator
$$
\eqalign{R_A\ :\, & \Hduale(\Omega)\longrightarrow\intersA{\Omega}\cr
		& \quad z\quad\,\longmapsto\qquad w\cr}
$$
where $w$ is the solution of
$$
\cases{Aw+\mu_A w=z\cr
	w\in\intersA{\Omega};\cr}\eqno(\frm{probleA})
$$
so that we have $u=R_Av$ (and $v=R_Bf$).

Define now the space
$$
V_A(\Omega):=\left\{w\in{\intersA{\Omega}}\,:\,\exists 
z\in{\intersA{\Omega}}\ \hbox{ s.t. } w=R_Az\right\}.
$$
Let us prove that the function $z$ is unique, for each $w$. 

Assume, by contradiction, that there exist two 
$z_1,\,z_2\in\intersA{\Omega}$ such that $w=R_Az_1=R_Az_2$. From the 
equation we have 
$$
\intl_\Omega (z_1-z_2)\psi\,dx=0,\qquad\forall\psi\in\intersA{\Omega}.
$$
Since $z_1-z_2\in\intersA{\Omega}$ we get 
$$
\intl_\Omega(z_1-z_2)^2\,dx=\Vert 
z_1-z_2\Vert^2_{\Ldue{\Omega}}=0
$$
and we conclude that $z_1=z_2$.

So we have defined
$$
V_A(\Omega)=R_A\left(\intersA{\Omega}\right)
$$
and on its image $R_A$ is one-to-one. So, on the space
$V_A(\Omega)$, the inverse operator can be defined:
$$
\eqalign{\Lambda_{\mu_C}\ :\, & 
V_A(\Omega)\longrightarrow\intersA{\Omega}.\cr
		& \quad w\quad\,\longmapsto\qquad z\ \hbox{ s.t. }w=R_Az\cr}
$$
Observe now that the operator $R_A$ is symmetric, that is:
$$
\dualita{h}{R_A g}\,=\,\dualita{g}{R_Ah},
$$
for any $h\in\Hduale(\Omega)$ and $g\in\Hunozero(\Omega)$.

We are now able to write system (\rf{problelimi}) as a single equation. We  
have $v=R_B f$ and $u=R_Av$, so that $u\in V_A(\Omega)$. For any $\varphi\in 
V_B(\Omega)$ 
there exists a unique $z\in\intersB{\Omega}$ such that $\varphi=R_Bz$. 
Hence, 
$$
\eqalign{
\dualita{f}{\varphi}&=\dualita{f}{R_Bz}\,=\,\dualita{z}{R_Bf}
	\,=\,\dualita{z}{v}\cr
	&=\dualita{\Lambda_{\mu_B}\varphi}{\Lambda_{\mu_A}u}
		=\intl_\Omega \Lambda_{\mu_B}\varphi\ \Lambda_{\mu_A}u\,dx.\cr}
$$
Last equality holds because both $\Lambda_{\mu_B}\varphi$ and 
$\Lambda_{\mu_A}u$ are in $\Ldue{\Omega}$.

Hence we can conclude that the function $u$, solution of (\rf{problelimi}), 
also solves
$$
\cases{
	u\in V_A(\Omega)\cr
	\displaystyle{\int_\Omega \Lambda_{\mu_B}\varphi\ \Lambda_{\mu_A}u\,dx}
		=\dualita{f}{\varphi}\cr
	\forall\varphi\in V_B(\Omega).\cr}
$$
It is easily seen that also the converse holds. Let $v$ be a solution of 
$$
\cases{Bv+\mu_B v=f\cr
	v\in\intersB{\Omega};\cr}
$$
then, by the definition, $v=R_Bf$. Hence, taken any $z\in\intersB{\Omega}$, 
we set $\varphi:=R_Bz$, and we have:
$$
\eqalign{\dualita{v}{z}=&\dualita{R_Bf}{z}\,
	=\,\dualita{f}{R_Bz}\,=\,\dualita{f}{\varphi}\cr
			=&\dualita{\Lambda_{\mu_A}u}{\Lambda_{\mu_B}\varphi}= 
				\dualita{\Lambda_{\mu_A}u}{z}.\cr}
$$
Taking then $z=v-\Lambda_{\mu_A}u$, we obtain
$$
\Vert v-\Lambda_{\mu_A}u\Vert^2_{{}_{\rm L^2}}=0,
$$
so 
$$
\Lambda_{\mu_A}u=v,
$$
that is, by the definition of $\Lambda_{\mu_A}$,
$$
u=R_Av,
$$
hence $u$ solves (\rf{problelimi}).

All this can be summarized in the following

\th{\thm{single}}{A function $u\in\intersA{\Omega}$ is a solution of
$$
\cases{
	u\in V_A(\Omega)\cr
	\displaystyle{\int_\Omega \Lambda_{\mu_B}\varphi\ \Lambda_{\mu_A}u\,dx}
		=\dualita{f}{\varphi}\cr
	\forall\varphi\in V_B(\Omega),\cr}
$$
if and only if it solves 
$$
\cases{Au+\mu_A u=v\cr
	u\in\intersA{\Omega}\cr
	Bv+\mu_B v=f\cr
	v\in\intersB{\Omega}.\cr}
$$}

It is importat to remark that the road back to the equivalence with a 
fourth order problem stops here. This is because, in 
general, the image of $\Lambda_{\mu_A}$ is not contained in the domain of 
$\Lambda_{\mu_B}$, even if $A=B$ and the two measures coincide.

\parag{\chp{optimization}}{The optimization problem}

Let us consider now the optimization problem (\rf{optim}).

Let $f$ be a function in $\Hduale(\Omega)$ and 
let $j:\Omega\times\R\rightarrow\R$ satisfy the standard Carath\'eodory 
conditions and be such that
$$
|j(x,s)|\leq b(x)+\beta |s|^p,\quad
\hbox{ for a.e. }x\in\Omega\hbox{ and }\forall s\in\R,\eqno(\frm{stima})
$$
with $b\in L^1(\Omega)$, $\beta\in\R$ and $1\leq p<2^*$, where
$2^*=\displaystyle{{2n}\over{n-2}}$ is the Sobolev exponent.

We want to study 
$$
\min_{U\in{\cal U}(\Omega)}\intl_\Omega j(x,u_{{}_U}(x))\,dx,\eqno(\frm{min})
$$
where ${\cal U}(\Omega)$ is the family of all open subsets of $\Omega$ and 
$u_{{}_U}$ is the solution, trivially extended in $\Omega\setminus U$, of the 
problem 
$$
\cases{\Delta^2u=f {\rm \ \ in\ } \Hduale(U)\cr
       \Delta u\in\Hunozero(U)\cr
       u\in\Hunozero(U).\cr}
$$
In general problem (\rf{min}) does not have solution, in the sense that
the infimum is not attained on the set
$\left\{u_{{}_U}\,:\,U\in{\cal 
U}(\Omega)\right\}$. This can be seen with the following example.

\ex{\thm{nosol}}{}Consider the case of $\Omega\subset\R^n$ with smooth 
boundary. Consider a function
$$
w\in C^\infty(\overline\Omega),\hbox{ s.t. } w>0\hbox{ in }\Omega,\ 
w=0\hbox{ and } \Delta w=0\hbox{ on }\partial\Omega  
$$
and set
$$
z:=-\Delta w+w,\quad f:=-\Delta z+z.\eqno(\frm{zf})
$$
Consider now 
$$
j(x,s)=\left(s-w(x)\right)^2,
$$
so that the integral
$$
\intl_\Omega \left(u(x)-w(x)\right)^2\,dx\eqno(\frm{integ})
$$
takes the value zero if and only if $u=w$ a.e. It will turn out from the
theory (see equation (\rf{infimum})) that zero is in fact the infimum of
the functional over all
solutionds
of Dirichlet problems on subsets of $\Omega$.

We prove now that, with this choice of $f$, $w$ can not be the solution 
of the problem
$$
\cases{\Delta^2 u_{{}_U}=f \hbox{ in }\Hduale(U)\cr
	\Delta u_{{}_U}\in\Hunozero(U)\cr
	u_{{}_U}\in\Hunozero(U),\cr}
$$
for any $U\subset\Omega$, or, 
equivalently (thanks to Proposition \rf{equiv}), of system
$$
\cases{-\Delta u_{{}_U}=v\hbox{ in }\Hduale(U)\cr
	u_{{}_U}\in\Hunozero(U)\cr
	-\Delta v=f\hbox{ in }\Hduale(U)\cr
	v\in\Hunozero(U).\cr}
$$ 

First consider the case when ${\rm cap}(\Omega\setminus U)=0$. Then the spaces 
$\Hunozero(\Omega)$ and $\Hunozero(U)$ coincide (see for instance Theorem 4.5 
in [\rf{HEI-KIL}]), and the problem in $U$ 
is the same as the problem in $\Omega$. We show that $w\neq u_\Omega$. 
From (\rf{zf}) we have
$$
\eqalign{\Delta^2w-\Delta w &=-\Delta z\cr
			&=f-z\cr
			&=f+\Delta w-w.\cr}
$$
So $\Delta^2w=f+2\Delta w -w$, and hence $w$ can be such that 
$\Delta^2w=f$ only if $-\Delta w=-\displaystyle{w\over 2}$ but this is impossible 
because, by the maximum principle, $w$ should be negative, while we have 
choosen it strictly positive.

The second case is when ${\rm cap}(\Omega\setminus U)>0$. If this happens 
$u_{{}_U}$ has to be zero 
in $\Omega\setminus U$, so that $u_{{}_U}\neq w$ on a set of non-zero 
capacity. But for functions in $\Hunozero(\Omega)$ to be equal
Lebesgue-a.e. is 
the same as cap-a.e., hence $u_{{}_U}$ can not be a minimizer for (\rf{integ}).

\medskip

The example shows then, that to be able to solve always our optimization 
problem, it is convenient to seek our minima in the larger set
$M:=\{u_\mu\,:\,\mu\in\Mo\}$, which is the closure of 
$N:=\{u_{{}_U}\,:\, U\in{\cal U}(\Omega)\}$ in the weak topology of
$\Hunozero(\Omega)$. In other words 
we introduce the relaxed optimization problem
$$
\min_{\mu\in\Mo}\intl_\Omega j(x,u_\mu(x))\,dx,\eqno(\frm{relaxopt})
$$
where $u_\mu$ denotes the solution fo the relaxed problem
$$
\cases{-\Delta u_\mu+u_\mu \mu=v_\mu\cr
	u_\mu\in\Hunozero(\Omega)\cap\Lduemu{\Omega}\cr
	-\Delta v_\mu+v_\mu \mu=f\cr
	v_\mu\in\Hunozero(\Omega)\cap\Lduemu{\Omega}\cr}\eqno(\frm{problemu})
$$ 
and each equation is meant in the sense of formula (\rf{relax}).

In order to see that solving the new problem (\rf{relaxopt}) gives the complete 
solution to problem (\rf{min}),
we have to show two things: that (\rf{relaxopt}) has a solution and that 
$$
\min_{\mu\in\Mo}\intl_\Omega j(x,u_\mu(x))\,dx\, 
	=\,\inf_{U\in{\cal U}(\Omega)}\intl_\Omega
j(x,u_{{}_U}(x))\,dx.\eqno(\frm{infimum})
$$

As for the first thing, we need to show that the set $M$ is closed in the 
weak topology of $\Hunozero(\Omega)$. So consider 
minimizing sequence $\mu_n$ and let $u_{\mu_n}$ be the solutions of the
systems
$$
\cases{-\Delta u_{\mu_n}+u_{\mu_n} \mu_n=v_{\mu_n}\cr
	u_{\mu_n}\in\Hunozero(\Omega)\cap L^2_{\mu_n}(\Omega)\cr
	-\Delta v_{\mu_n}+v_{\mu_n} \mu_n=f\cr
	v_{\mu_n}\in\Hunozero(\Omega)\cap L^2_{\mu_n}(\Omega).\cr}
$$ 

Applying Theorem \rf{compactness} to the second equation we have that
there exist a subsequence, called for simplicity 
$\mu_n$, and a measure $\mu\in\Mo$ such that $v_{\mu_n}\wto v_\mu$ weakly 
in $\Hunozero(\Omega)$, 
$$
\cases{-\Delta v_\mu+v_\mu \mu=f\cr
	v_\mu\in\Hunozero(\Omega)\cap\Lduemu{\Omega}.\cr}
$$

Applying again Theorem \rf{compactness} to the first equation (only
those corresponding to the subsequence), we obtain that 
$$
u_{\mu_n}\wto u_\mu\hbox{ weakly in }\Hunozero(\Omega).
$$
and $u_\mu$ is the solution of (\rf{problemu}).

Hence $u_\mu$  is admissible and, by the minimizing property of the sequence 
$\{\mu_n\}$,
$$
\intl_\Omega j(x,u_\mu(x))\,dx=\min_{\nu\in\Mo}\intl_\Omega 
j(x,u_\nu(x))\,dx.
$$

On the other hand observe that every problem of the form 
$$
\cases{-\Delta z=w\cr
	z\in\Hunozero(U)\cr
	-\Delta w=g\cr
	w\in\Hunozero(U)\cr}
$$
can be viewed as a relaxed problem with the measure
$$
\mu^U(B):=\cases{0\quad&if $B\setminus U$ has capacity zero\cr
		+\infty&otherwise\cr}
$$
since this holds for second order problems as can be seen, for 
instance, in [\rf{BUT-DM}].

This ensures us that problem (\rf{relaxopt}) is indeed an extension of 
(\rf{min}), in the sense that $N\subseteq M$, so that
$$
\min_{\mu\in\Mo}\intl_\Omega j(x,u_\mu(x))\,dx\leq
	\inf_{U\in{\cal U}(\Omega)}\intl_\Omega j(x,u_{{}_U}(x))\, dx.
$$

To prove the inverse inequality we need the following 
\lemma{\thm{density4}}{Let $\mu\in\Mo$. Then there exists a sequence $U_n$ 
of subsets of $\Omega$ such that
$$
u_{{}_{U_n}}\wto u_\mu\hbox{ weakly in }\Hunozero(\Omega)
$$
where $u_{{}_{U_n}}$ and $u_\mu$ solve respectively
$$
\cases{-\Delta u_{{}_{U_n}}=v_n\hbox{ in }\Hduale(U_n)\cr
	u_{{}_{U_n}}\in\Hunozero(U_n)\cr
	-\Delta v_n=f\hbox{ in }\Hduale(U_n)\cr
	v_n\in\Hunozero(U_n)\cr}\qquad\hbox{ and }
\cases{-\Delta u_\mu+u_\mu \mu=v_\mu\cr
	u_\mu\in\Hunozero(\Omega)\cap\Lduemu{\Omega}\cr
	-\Delta v_\mu+v_\mu \mu=f\cr
	v_\mu\in\Hunozero(\Omega)\cap\Lduemu{\Omega}.\cr}
$$ }
\proof
We know, from Theorem \rf{density}, that given $\mu$, we have a 
sequence of sets $U_n$ such that, for every $g\in\Hduale(\Omega)$, 
if 
$$
\cases{-\Delta z_n=g\hbox{ in }\Hduale(U_n)\cr
	z_n\in\Hunozero(U_n).\cr}
$$
then the sequence $z_n$ tends weakly in $\Hunozero(\Omega)$  to $z_\mu$, 
solution of 
$$
\cases{ -\Delta z_\mu+v_\mu \mu=g\cr
	z_\mu\in\Hunozero(\Omega)\cap\Lduemu{\Omega}.\cr}
$$
Since the sequence $U_n$ depends only on the operator and not on the right hand 
side the same sequence of sets suits for the problem in $u$. Hence 
$u_{{}_{U_n}}$ converge to $u_\mu$ solution of the second system.
\endproof

Since, by the hypothesis (\rf{stima}), the operator
$u\mapsto{\displaystyle \int 
j(x,u)}\,dx$ is continuous in the strong topology of ${\rm L}^2(\Omega)$
(see 
[\rf{BUT-DM}] Theorem 4.1), the density result enables us to conclude that
$$
\min_{\mu\in\Mo}\intl_\Omega j(x,u_\mu(x))\,dx\,=\,
	\inf_{U\in{\cal U}(\Omega)}\intl_\Omega j(x,u_A(x))\, dx.
$$
This ensures us also that every solution of a relaxed system can be 
approximated by solutions of fourth order problems.

\rem{\thm{estensione}}{}We remark here that the model case of the 
bi-laplacian can be extended without changes in the proofs to the case of 
a fourth order elliptic operator with constat coefficients $L$ such that
$$
Lu\,=\,C^2u,
$$
where $C$ is a second order elliptic operator with constant coefficients.

\intro{References}
\def\interrefspace{\smallskip}  
{\ninepoint\frenchspacing

\item{[\bib{ATT}]}ATTOUCH H.: {\it Variational converence for functions
and operators}. Pitman, Boston, 1984 
\interrefspace

\item{[\bib{BUT-DM}]}BUTTAZZO  G., DAL MASO  G.: Shape optimization 
for Dirichlet Problems: Relaxed Formulation and Optimality Conditions. 
{\it Appl. Math. Optim. \/} {\bf 23} (1991), 17-49. 
\interrefspace

\item{[\bib{CHI-DM}]}CHIPOT  M., DAL MASO  G.: Relaxed shape 
optimization: thecase of nonnegative data for the Dirichlet 
problem. {\it Adv. Math. Sci. Appl.\/} {\bf 1} (1992), 47-81.
\interrefspace

\item{[\bib{CIA1}]}CIARLET P.G.: {\it The finite element method for 
elliptic problems}. Studies in mathematics and its applications. North 
Holland, New York, 1978 
\interrefspace

\item{[\bib{CIA2}]}CIARLET P.G.: {\it Plates and junctions in elastic 
multi-structures. An asymptotic analysis}. Masson, paris, 1990 \interrefspace

\item{[\bib{CIO-MUR}]}CIORANESCU  D., MURAT  F.: Un terme \'etrange venu 
d'ailleurs I and II. In {\it Nonlinear partial differential equations and 
their applications.} Pitman research notes in mathematics; n. 60, 98-138; 
n. 70, 154-178. 
\interrefspace

\item{[\bib{DM-GAR}]}DAL MASO  G., GARRONI  A.: New results on 
the asymptotic behaviour of Dirichlet problems in perforated domains.
{\it Math. Mod. Meth. Appl. Sci.\/} {\bf 3} (1994), 373-407.
\interrefspace

\item{[\bib{DM-MAL}]}DAL MASO  G., MALUSA  A.: Approximation of relaxed 
Dirichlet problems by boundary value problems in perforated domains. {\it 
Proc. Roy. Soc. Edinburgh\/} {\bf 125A} (1995), 99-114.
\interrefspace

\item{[\bib{DM-MOS}]}DAL MASO  G., MOSCO  U.: Wiener's criterion and 
$\Gamma$-convergence. {\it Appl. Math. Optim.\/} {\bf 15} (1987),
15-63.
\interrefspace

\item{[\bib{DM-PAD}]}DAL MASO  G., PADERNI  G.: Variational inequalities 
for the biharmonic operator with variable coefficients. {\it Ann. Mat. 
Pura Appl.\/} {\bf 103} (1988), 203-227. 
\interrefspace

\item{[\bib{DUV-LIO}]}DUVAUT  G., LIONS  J.-L.: {\it Inequalities in 
mechanics and physics.}  Springer-Verlag, Berlin, 1976
\interrefspace

\item{[\bib{EVA-GAR}]}EVANS  L.C., GARIEPY  R.F.: {\it Measure theory and 
fine properties of functions.} CRC Press, Boca Raton, 1992
\interrefspace

\item{[\bib{HEI-KIL}]}HEINONEN  J., KILPEL\"AINEN  T., MARTIO  O.: {\it 
Nonlinear potential theory of degenerate elliptic equations.} Clarendon 
Press, Oxford, 1993

\interrefspace

\item{[\bib{LAG-LIO}]}LAGNESE  J.E., LIONS  J.-L.: {\it Modelling 
analysis and control of thin plates.} Masson, Paris, 1988 \interrefspace

\item{[\bib{LAN-LIF}]}LANDAU  L., LIFSHITZ  E.: {\it Theory of 
elasticity.} Course of theoretical physics. Pergamon Press, 1975 
\interrefspace

\end